\documentclass[11pt,a4paper,thmdefs]{amsart}
\usepackage{amsfonts}
\usepackage{amssymb,amsmath}
\usepackage{amscd}
\usepackage{graphicx}

\newtheorem{theorem}{Theorem}[section]
\newtheorem{lemma}{Lemma}[section]

\newtheorem{corollary}{Corollary}[section]
\newtheorem{definition}{Definition}[section]
\newtheorem{example}{Example}[section]

\numberwithin{equation}{section} \textwidth12.5cm
\input{tcilatex}

\begin{document}
\title[]{On a new type Bertrand curve}
\author[]{\c{C}ET\.{I}N CAMCI}
\address{Department of Mathematics \\
Onsekiz Mart University \\
17020 \c{C}anakkale, Turkey}
\email{ccamci@comu.edu.tr}
\date{}

\begin{abstract}
In this paper, we define a new type Bertrand curve and this curves are said $%
V$-Bertrand curve, $f$-Bertrand curve and $a$-Bertrand curve. In addition,
we give charectarization of the $V$-Bertrand curve and we define a Bertrand
surface.
\end{abstract}

\maketitle


\section{\textbf{Introduction}}

Let $\gamma :I\longrightarrow 
%TCIMACRO{\U{211d} }%
%BeginExpansion
\mathbb{R}
%EndExpansion
^{3}$ $\left( s\rightarrow \gamma (s)\right) $ be a unit-speed curve.
Formula of the Serret-Frenet equations are given by \ 
\begin{equation*}
\left( 
\begin{array}{c}
T%
%TCIMACRO{\U{b4} }%
%BeginExpansion
{\acute{}}
%EndExpansion
\\ 
N^{%
%TCIMACRO{\U{b4}}%
%BeginExpansion
{\acute{}}%
%EndExpansion
} \\ 
B^{%
%TCIMACRO{\U{b4}}%
%BeginExpansion
{\acute{}}%
%EndExpansion
}%
\end{array}%
\right) =\left( 
\begin{array}{ccc}
0 & \kappa & 0 \\ 
-\kappa & 0 & \tau \\ 
0 & -\tau & 0%
\end{array}%
\right) \left( 
\begin{array}{c}
T \\ 
N \\ 
B%
\end{array}%
\right)
\end{equation*}%
where $\left\{ T,N,B,\kappa ,\ \tau \right\} $ \ is the Serret-Frenet
apparatus of the curves(\cite{F}, \cite{S}) .

In 3-Euclidean spaces, Saint-Venant proposed that the ruled surface whose
based curve is the principal normals of a curve and another curve which have
same principal normals can exist (\cite{SV}). Solution of the problem was
given by Bertrand(\cite{JB}). Let $M$ be a unit-speed curve with coordinate
neighborhood $(I,\gamma )$ of the curve. So we can define a curve by \ 
\begin{equation*}
\beta (s)=\gamma (s)+\lambda (s)N(s)
\end{equation*}%
where $\lambda :I\longrightarrow 
%TCIMACRO{\U{211d} }%
%BeginExpansion
\mathbb{R}
%EndExpansion
$ $\left( s\rightarrow \lambda (s)\right) $ is differentiable function. Let $%
\left\{ \overline{T},\overline{N},\overline{B}\right\} $ be Frenet frame of
the curve and $\overline{\kappa }$ and \ $\overline{\tau }$ be curvatures of
the curves $K$. \ If $\left\{ N,\overline{N}\right\} $ is linear dependent $%
\left( \overline{N}=\epsilon N,\epsilon =\pm 1\right) $, then it is said
that $M$ is Bertrand mate of $K$(\cite{DC},\cite{HCC}). \ If there is
Bertrand mate of the curve $\ M$, then it is said that $M$ is Bertrand curve
(\cite{DC},\cite{HCC}). For example, any planar curve and circular helix is
Bertrand curve. It is well known theorem that $M$ is Bertrand curve if and
only if there exist $\lambda $, $\mu \in 
%TCIMACRO{\U{211d} }%
%BeginExpansion
\mathbb{R}
%EndExpansion
$ such that $\lambda \kappa +\mu \tau =1$ where $\lambda \neq 0$ , $\mu
=\lambda \cot \theta $ and $\theta $ is angle between $T$ and $\overline{T}$
(\cite{DC},\cite{HCC}). In fact, $M$ is Bertrand or $B-$Bertrand curve if
and only if there are $\lambda $, $\mu \in 
%TCIMACRO{\U{211d} }%
%BeginExpansion
\mathbb{R}
%EndExpansion
$ such that $\lambda \kappa +\mu \tau =1$.\ 

If the curve is a unit-speed spherical curve, for all $s\in I$, $\gamma (s)$
is perpendicular to $T=\gamma ^{\prime }(s)$. So we have Sabban frame as 
\begin{equation*}
\left\{ \gamma (s)\text{, }T(s)=\gamma ^{\prime }(s)\text{, }Y(s)=\gamma
(s)\times T(s)\right\} 
\end{equation*}%
(\cite{ITB}). Serret-Frenet formula of the this curve is given by 
\begin{equation*}
\left( 
\begin{array}{c}
\gamma ^{\prime } \\ 
T^{\prime } \\ 
Y^{\prime }%
\end{array}%
\right) =\left( 
\begin{array}{ccc}
0 & 1 & 0 \\ 
-1 & 0 & \kappa _{g} \\ 
0 & -\kappa _{g} & 0%
\end{array}%
\right) \left( 
\begin{array}{c}
\gamma  \\ 
T \\ 
Y%
\end{array}%
\right) 
\end{equation*}%
where $\kappa _{g}(s)=\det \left( \gamma ,T,T^{\prime }\right) $ (\cite{ITB}%
). Under the above notation, the space curve is definned by 
\begin{equation}
\widetilde{\gamma }(s)=a\int \gamma (s)ds+a\cot \theta \int Y(s)ds+c
\label{x}
\end{equation}%
where $a,\theta \in 
%TCIMACRO{\U{211d} }%
%BeginExpansion
\mathbb{R}
%EndExpansion
$ and $c$ is constant vector (\cite{ITB}). Izumiya and Tkeuchi prove that
all Bertrand curves can be constructed from equation (\ref{x}) and they gave
some Bertrand curve examples by using equation (\ref{x}) (\cite{ITB}). In
this paper, we give two method different from above method. Furthermore, we
define a $V$-Bertrand curve, $f$-Bertrand curve and $a$-Bertrand curve where
Bertrand curve is $T-$Bertrand curve and $a\in 
%TCIMACRO{\U{211d} }%
%BeginExpansion
\mathbb{R}
%EndExpansion
$. Furthermore, we modify charectarization of the Bertrand curve. In
addition, we define a Bertrand surface and we give an example Bertrand
surface. Moreover, we define a equivalence relation on Bertrand surface set.

\section{\textbf{Preliminaries}}

In 3-Euclidean spaces, let $M$ be a unit-speed curve and $\left\{
T,N,B,\kappa ,\ \tau \right\} $ \ is Serret-Frenet apparatus of the curves.
The curve lies on 2-sphere if and only if below equation is hold 
\begin{equation}
\left( \left( \frac{1}{\kappa }\right) ^{\prime }\frac{1}{\tau }\right)
^{\prime }+\frac{\tau }{\kappa }=0  \label{1.1}
\end{equation}%
(\cite{EK},\cite{ST}). From the solution of the equation (\ref{1.1}), we
have 
\begin{equation}
\frac{1}{\kappa }=R\cos \left( \overset{s}{\underset{0}{\dint }}\tau
(u)du+\theta _{0}\right) .  \label{1.4}
\end{equation}%
where radius of the sphere is equal to $R$ ( \cite{EK},\cite{BG}, \cite{wong}%
). It is well known that if $M$ is spherical curve, then $M$ and osculating
circle lies on this sphere. From equation (\ref{1.4}), we have 
\begin{equation*}
\frac{1}{\kappa _{0}}=\cos \theta _{0}=\frac{R_{0}}{R}.
\end{equation*}%
If $R$ is equal to $1$, then we have $\cos \theta _{0}=R_{0}=\frac{1}{\kappa
_{0}}$ (\cite{CC}).

Choi and Kim defined a unit vector field $V$ which is given by 
\begin{equation*}
V(s)=u(s)T(s)+v(s)N(s)+w(s)B(s)
\end{equation*}%
and integral curve of $V$ is define by%
\begin{equation*}
\gamma _{V}(s)=\int V(s)ds
\end{equation*}%
where $u,v,w$ are functions from $I$ to $%
%TCIMACRO{\U{211d} }%
%BeginExpansion
\mathbb{R}
%EndExpansion
$ and $u^{2}+v^{2}+w^{2}=1$ (\cite{Choi}). Without loss of generality, we
suppose that arc-lenght parameter of $M$ and integral curve of $V$ are same (%
\cite{Choi}). If $u(s)=1,v(s)=w(s)=0$, for all $s\in I$, then we have 
\begin{equation*}
\gamma _{T}(s)=\int T(s)ds=\gamma (s).
\end{equation*}%
If $u(s)=w(s)=0$, $v(s)=1$, for all $s\in I$, then we have 
\begin{equation*}
\gamma _{N}(s)=\int N(s)ds
\end{equation*}%
and if $u(s)=v(s)=0$, $w(s)=1$, for all $s\in I$, then we have 
\begin{equation*}
\gamma _{B}(s)=\int B(s)ds
\end{equation*}%
(\cite{Choi}). $\left( \gamma _{N}\right) $ (resp.$\left( \gamma _{B}\right) 
$) is called principal-direction curve (resp. binormal-direction curve) of $%
\gamma $ (\cite{Choi}). Let 
\begin{equation*}
\kappa _{1}=\sqrt{\kappa ^{2}+\tau ^{2}},\tau _{1}=\frac{\kappa ^{2}}{\kappa
^{2}+\tau ^{2}}\left( \frac{\tau }{\kappa }\right) ^{\prime }
\end{equation*}%
and $T_{1}=N$, $N_{1}=\frac{N^{\prime }}{\left\Vert N^{\prime }\right\Vert }$%
, $B_{1}=T_{1}\times N_{1}$ ( resp. $\kappa _{\beta }=\left\vert \tau
\right\vert ,\tau _{\beta }=\kappa $ and $T_{\beta }=B,N_{\beta }=N,B_{\beta
}=-T$ ) be curvatures and Frenet vectors of principal-direction curve (resp.
binormal-direction curve)(\cite{Choi}, \cite{Mont}). If $\ $%
\begin{equation*}
u(s)=0,v(s)=-\cos \left( \int \tau (s)ds\right) \neq 0,w(s)=\sin \left( \int
\tau (s)ds\right)
\end{equation*}%
then $\gamma _{pdo}(s)=\int V(s)ds$ is called principal-donor curve of $%
\gamma $, for all $s\in I$ (\cite{Choi}). So we have 
\begin{equation*}
\kappa _{pdo}=\kappa \cos \left( \int \tau (s)ds\right) ,\tau _{pdo}=\kappa
\sin \left( \int \tau (s)ds\right)
\end{equation*}%
and $T_{pdo}=V$, $N_{pdo}=T$, $B_{pdo}=T_{pdo}\times N_{pdo}$ .

\section{V-Bertrand Curve }

Let $\gamma :I\longrightarrow 
%TCIMACRO{\U{211d} }%
%BeginExpansion
\mathbb{R}
%EndExpansion
^{3}$ $\left( s\rightarrow \gamma (s)\right) $ be a unit-speed curve and $%
\left\{ T,N,B,\kappa ,\ \tau \right\} $ \ be Serret-Frenet apparatus of the
curve. We can define a curve $K$ by 
\begin{equation*}
\beta (s)=\int V(s)ds+\lambda (s)N(s)
\end{equation*}%
where $\lambda :I\longrightarrow 
%TCIMACRO{\U{211d} }%
%BeginExpansion
\mathbb{R}
%EndExpansion
$ $\left( s\rightarrow \lambda (s)\right) $ is differentiable function and $%
V=uT+vN+wB$ is unit vector field. Let $\left\{ \overline{T},\overline{N},%
\overline{B}\right\} $ be orthonormal frame of the curve and $\overline{%
\kappa }$ and \ $\overline{\tau }$ be curvatures of the curve $K$. So we can
give following definition.

\begin{definition}
If $\left\{ N,\overline{N}\right\} $ is linear dependent $\left( \overline{N}%
=\epsilon N,\epsilon =\pm 1\right) $, then it is said that $\left(
M,K\right) $ is $V$-Bertrand mate. If $V=T$, then $\left( M,K\right) $ is
Bertrand mate.

\begin{theorem}
Let $M$ be a unit-speed curve and $\left\{ T,N,B,\kappa ,\ \tau \right\} $ \
be Serret-Frenet apparatus of the curves. $M$ is $V$-Bertrand curve if and
only if it is satisfied that 
\begin{equation}
\lambda \left( \kappa \tan \theta +\tau \right) =u\tan \theta -w  \label{3.1}
\end{equation}%
and 
\begin{equation}
\lambda (s)=-\int v(s)ds  \label{3.2}
\end{equation}%
where $\theta $ is constant angle between $T$ and $\overline{T}$.

\begin{proof}
If $M$ is $V$-Bertrand curve, then $V$-Bertrand mate of $M$ is equal to 
\begin{equation}
\beta (s)=\int V(s)ds+\lambda (s)N(s)  \label{3.3}
\end{equation}%
and $\left\{ N,\overline{N}\right\} $ is linear dependent. \ If we
derivative equation (\ref{3.3}), we have 
\begin{equation*}
\frac{d\overline{s}}{ds}\overline{T}=\left( u-\lambda \kappa \right)
T+\left( \lambda ^{\prime }+v\right) N+\left( w+\lambda \tau \right) B.
\end{equation*}%
Since $\left\{ N,\overline{N}\right\} $ is linear dependent, we have 
\begin{equation*}
\lambda (s)=-\int v(s)ds
\end{equation*}%
and we get 
\begin{equation}
\overline{T}=\frac{ds}{d\overline{s}}\left( u-\lambda \kappa \right) T+\frac{%
ds}{d\overline{s}}\left( w+\lambda \tau \right) B=\cos \theta (s)T+\sin
\theta (s)B  \label{3.4}
\end{equation}%
where $\cos \theta (s)=\frac{ds}{d\overline{s}}\left( u-\lambda \kappa
\right) $ and $\sin \theta (s)=\frac{ds}{d\overline{s}}\left( w+\lambda \tau
\right) $. So we have 
\begin{equation*}
\tan \theta =\frac{w+\lambda \tau }{u-\lambda \kappa }
\end{equation*}%
or 
\begin{equation*}
\lambda \left( \kappa \tan \theta +\tau \right) =u\tan \theta -w\text{.}
\end{equation*}%
If we take the derivative of a equation (\ref{3.4}), we have 
\begin{equation*}
\frac{d\overline{s}}{ds}\overline{\kappa }\overline{N}=-\theta ^{\prime
}\sin \theta T+\left( \kappa \cos \theta -\tau \sin \theta \right) N+\theta
^{\prime }\cos \theta B
\end{equation*}%
So that $\left\{ N,\overline{N}\right\} $ is linear dependent, we have $%
\theta ^{\prime }=0$ or $\theta =const$. Conversely, we define a curve by 
\begin{equation}
\beta (s)=\int V(s)ds+\lambda (s)N(s)  \label{3.5}
\end{equation}%
where $\lambda :I\longrightarrow 
%TCIMACRO{\U{211d} }%
%BeginExpansion
\mathbb{R}
%EndExpansion
$ $\left( s\rightarrow \lambda (s)\right) $ is differentiable function such
that $\lambda (s)=-\int v(s)ds$ . If we derivative equation (\ref{3.5}), we
have%
\begin{equation*}
\frac{d\overline{s}}{ds}\overline{T}=\left( u-\lambda \kappa \right)
T+\left( \lambda ^{\prime }+v\right) N+\left( w+\lambda \tau \right) B.
\end{equation*}%
From equation (\ref{3.2}), we get 
\begin{equation}
\overline{T}=\frac{ds}{d\overline{s}}\left( u-\lambda \kappa \right) T+\frac{%
ds}{d\overline{s}}\left( w+\lambda \tau \right) B=\cos \phi (s)T+\sin \phi
(s)B  \label{3.6}
\end{equation}%
where $\cos \phi (s)=\frac{ds}{d\overline{s}}\left( u-\lambda \kappa \right) 
$ and $\sin \phi (s)=\frac{ds}{d\overline{s}}\left( v+\lambda \tau \right) $%
. From equation (\ref{3.1}) and (\ref{3.6}) we have 
\begin{equation*}
\tan \phi (s)=\frac{v+\lambda \tau }{u-\lambda \kappa }=\tan \theta .
\end{equation*}%
\ 

So $\phi =\theta +n\pi $ is angle between $T$ and $\overline{T}$ where $n\in
Z$. If we take the derivative of a equation (\ref{3.6}), we see that $%
\left\{ N,\overline{N}\right\} $ is linear dependent.
\end{proof}

\begin{definition}
Let $M$ be a unit-speed non-planar curve and $\left\{ T,N,B,\kappa ,\ \tau
\right\} $ \ is Serret-Frenet apparatus of the curves. If there is $\lambda
\neq 0$ , $\theta \in 
%TCIMACRO{\U{211d} }%
%BeginExpansion
\mathbb{R}
%EndExpansion
$ such that 
\begin{equation*}
\lambda \kappa +\lambda \cot \theta \tau =1,
\end{equation*}%
then we said that $M$ is Bertrand curve (or $T$-Bertrand curve ). If there
is $\lambda \neq 0$ , $\theta \in 
%TCIMACRO{\U{211d} }%
%BeginExpansion
\mathbb{R}
%EndExpansion
$ such that 
\begin{equation*}
\lambda \tan \theta \kappa +\lambda \tau =-1,
\end{equation*}%
then we said that $M$ is $B$-Bertrand curve.
\end{definition}
\end{theorem}

\begin{corollary}
If $u(s)=1,v(s)=w(s)=0$, we have Bertrand curve. From equation (\ref{3.2}), $%
\lambda $ is constant. \ From equation (\ref{3.1}), we have $\lambda \kappa
+\mu \tau =1$ where $\lambda \neq 0$ , $\mu =\lambda \cot \theta $ \ and $\
\theta $\ is constant angle between $T$ and $\overline{T}$. If $%
v(s)=1,u(s)=w(s)=0$, we have $B$-Bertrand curve. From equation (\ref{3.1}),
we have $\lambda \kappa +\mu \tau =-1$ where $\mu \neq 0$ is constant, $%
\lambda =\mu \tan \theta $ and $\theta $ is constant angle between $T$ and $%
\overline{T}$.
\end{corollary}

\begin{corollary}
If the curve is a Salkowki curve, then curvature of the curve is equal to $1$%
(or constant). If $\lambda =1$ , $\mu =\lambda \cot \theta =0$, then we have 
$\lambda \kappa +\mu \tau =1$. \ if $\mu =0$ , $\lambda =\mu \cot \theta =-1$%
, then we have $\lambda \kappa +\mu \tau =-1$. \ But this is contradiction.
So a Salkowski curve is a Bertrand curve, but it is not $B$-Bertrand curve.
If the curve is anti-Salkowki curve, then torsion of the curve is equal to $1
$(or constant). If $\mu =-1$ , $\lambda =\mu \tan \theta =0$, then we have $%
\lambda \kappa +\mu \tau =-1$. \ if $\lambda =0$ , $\mu =\lambda \cot \theta
=-1$, then we have $\lambda \kappa +\mu \tau =-1$. \ But this is a
contradiction. So a anti-Salkowski curve is a $B$-Bertrand curve, but it is
not a Bertrand curve.
\end{corollary}

\begin{corollary}
If the curve is planar curve ($\tau =0$), then the curve is not only
Bertrand but also $B$-Bertrand curve.
\end{corollary}

\begin{corollary}
If $v(s)=1,u(s)=w(s)=0$, we have $N$-Bertrand curve. From equation (\ref{3.2}%
),we have $\lambda =c-s$ . From equation (\ref{3.1}), we have $\frac{\tau }{%
\kappa }=-\tan \theta $ where $\theta $ is constant angle between $T$ and $%
\overline{T}$. As a result, $M$ is $N$-Bertrand curve if and only if $M$ is
general helix.
\end{corollary}

\begin{corollary}
If the curve $M$ is Bertrand curve, there is $\lambda \neq 0$, $\mu \neq
0\in 
%TCIMACRO{\U{211d} }%
%BeginExpansion
\mathbb{R}
%EndExpansion
$ \ such that $\lambda \kappa +\mu \tau =1$ . If we define $\lambda
_{1}=-\lambda $ and $\mu _{1}=$ $-\mu $, then we have $\lambda _{1}\kappa
+\mu _{1}\tau =-1$. As a result, the curve $M$ is $B$-Bertrand curve.
Conversely, similarly, If the curve $M$ is $B$-Bertrand curve where there
exist $\lambda \neq 0$, $\mu \neq 0\in 
%TCIMACRO{\U{211d} }%
%BeginExpansion
\mathbb{R}
%EndExpansion
$ \ such that $\lambda \kappa +\mu \tau =-1$ , then the curve $M$ is
Bertrand curve. In this case, it is said that $M$ is proper Bertrand curve
\end{corollary}

\begin{corollary}
Let $\left( \alpha \right) $ be a unit-speed curve and $\left\{ T,N,B,\kappa
,\ \tau \right\} $ \ be Serret-Frenet apparatus of the curves. We suppose
that $v=0$, then we can see that $\lambda $ is non zero constant$.$ So we
get $u^{2}+w^{2}=1$ and $w=\varepsilon \sqrt{1-u^{2}}$where $\varepsilon
=\pm 1$. From equation (\ref{3.1}), we have%
\begin{equation}
u\tan \theta -\varepsilon \sqrt{1-u^{2}}=f  \label{3.7}
\end{equation}%
So, from (\ref{3.7}), we have 
\begin{equation}
u^{\pm }=\frac{f\tan \theta \pm \sqrt{1+\left( \tan \theta \right) ^{2}-f^{2}%
}}{1+\left( \tan \theta \right) ^{2}}\text{.}  \label{3.8}
\end{equation}%
If $u^{+}=\frac{f\tan \theta +\sqrt{1+\left( \tan \theta \right) ^{2}-f^{2}}%
}{1+\left( \tan \theta \right) ^{2}}$, then we have 
\begin{equation*}
w_{1}^{\pm }=\pm \sqrt{1-\left( u^{+}\right) ^{2}}
\end{equation*}%
If $u^{-}=\frac{f\tan \theta -\sqrt{1+\left( \tan \theta \right) ^{2}-f^{2}}%
}{1+\left( \tan \theta \right) ^{2}}$, then we have 
\begin{equation*}
w_{2}^{\pm }=\pm \sqrt{1-\left( u^{-}\right) ^{2}}
\end{equation*}%
So, we have 
\begin{eqnarray*}
V_{1}^{+} &=&u^{+}T+w_{1}^{+}B,V_{1}^{-}=u^{+}T+w_{1}^{-}B, \\
V_{2}^{+} &=&u^{-}T+w_{2}^{+}B,V_{2}^{-}=u^{-}T+w_{2}^{-}B
\end{eqnarray*}%
and 
\begin{eqnarray*}
\beta _{1}^{+} &=&\int V_{1}^{+}ds+\lambda N,\beta _{1}^{-}=\int
V_{1}^{-}ds+\lambda N \\
\beta _{2}^{+} &=&\int V_{2}^{+}ds+\lambda N,\beta _{2}^{-}=\int
V_{2}^{-}ds+\lambda N
\end{eqnarray*}%
In this case, $\left( \alpha \right) $ is $V_{1}^{+},V_{1}^{-},V_{2}^{+}$
and $V_{2}^{-}-$Bertrand curve. So we can give following definition.
\end{corollary}
\end{definition}

\begin{definition}
The curves $\beta _{1}^{+},\beta _{1}^{-},\beta _{2}^{+},\beta _{2}^{-}$ are
called $f-$Bertrand curve mate of $\left( \alpha \right) $ and the curve $%
\left( \alpha \right) $ is called $f-$Bertrand curve. If $f=r(const.)$, \
then $\left( \alpha \right) $ is called $r-$Bertrand curve and the curves $%
\beta _{1}^{+},\beta _{1}^{-},\beta _{2}^{+},\beta _{2}^{-}$ are called $r-$%
Bertrand curve mate of $\left( \alpha \right) $

\begin{example}
\ If \ $\alpha (s)=(a\cos \frac{s}{c},a\sin \frac{s}{c},\frac{bs}{c})$, then
we have 
\begin{eqnarray*}
T &=&(-\frac{a}{c}\sin \frac{s}{c},\frac{a}{c}\cos \frac{s}{c},\frac{b}{c})
\\
N &=&(-\cos \frac{s}{c},-\sin \frac{s}{c},0) \\
B &=&(\frac{b}{c}\sin \frac{s}{c},-\frac{b}{c}\cos \frac{s}{c},\frac{a}{c})
\end{eqnarray*}%
and $\kappa =\frac{a}{c^{2}}$and $\tau =\frac{b}{c^{2}}$ where $%
c^{2}=a^{2}+b^{2}$. If $\overline{\lambda }=\lambda \tan \theta $ and $%
\overline{\mu }=\lambda $, then we have $\overline{\lambda }\kappa +%
\overline{\mu }\tau =1=f$ where $\tan \theta =\frac{c^{2}-\lambda b}{\lambda
a}$. From (\ref{3.8}), we have 
\begin{equation*}
u^{\pm }=\frac{\tan \theta \pm \tan \theta }{1+\left( \tan \theta \right)
^{2}}
\end{equation*}%
So we have $u^{-}=0$ and $u^{+}=\frac{2\lambda a\left( c^{2}-\lambda
b\right) }{c^{2}\left( a^{2}+\left( \lambda -b\right) ^{2}\right) }$ . If $%
u^{-}=0$, then we have $w_{2}^{+}=1$ and $\ w_{2}^{-}=-1$. If $u^{+}=\frac{%
2\lambda a\left( c^{2}-\lambda b\right) }{c^{2}\left( a^{2}+\left( \lambda
-b\right) ^{2}\right) }$, then we have 
\begin{equation*}
w_{1}^{+}=\frac{\lambda ^{2}a^{2}+c^{2}\left( 2\lambda b-c^{2}\right) }{%
c^{2}\left( a^{2}+\left( \lambda -b\right) ^{2}\right) }
\end{equation*}%
and $\ $%
\begin{equation*}
w_{1}^{-}=-\frac{\lambda ^{2}a^{2}+c^{2}\left( 2\lambda b-c^{2}\right) }{%
c^{2}\left( a^{2}+\left( \lambda -b\right) ^{2}\right) }
\end{equation*}%
So we have 
\begin{equation*}
V_{2}^{+}=B,V_{2}^{-}=-B,V_{1}^{+}=u^{+}T+w_{2}^{+}B,V_{1}^{-}=u^{-}T+w_{2}^{-}B
\end{equation*}%
and 
\begin{eqnarray*}
\beta _{2}^{\pm } &=&\left( (\lambda \mp b)\cos \frac{s}{c},(\lambda \mp
b)\sin \frac{s}{c},\pm \frac{as}{c}\right)  \\
\beta _{1}^{\pm } &=&\left( \left( w_{2}^{\pm }b-u_{2}a-\lambda \right) \cos 
\frac{s}{c},\left( w_{2}^{\pm }b-u_{2}a-\lambda \right) \sin \frac{s}{c},%
\frac{w_{2}^{\pm }a+u_{2}b}{c}s\right) .
\end{eqnarray*}%
So, $\left( \alpha \right) $ is $1-$Bertrand curve and $\beta _{1}^{+},\beta
_{1}^{-},\beta _{2}^{+},\beta _{2}^{-}$ curves are all $1-$Bertrand curve
mate of $\left( \alpha \right) $.
\end{example}

\begin{example}
\ From (\ref{3.8}), if \ $\overline{\lambda }\kappa +\overline{\mu }\tau
=\tan \theta =f$ , then we obtain $u^{+}=1$ and $u^{-}=-\cos 2\theta $ where 
$\tan \theta =\frac{c^{2}-\lambda b}{\lambda a}$. If $u^{+}=1$ (resp.$%
u^{-}=-\cos 2\theta $), then we have $w_{2}^{\pm }=0$ (resp. $w_{1}^{\pm
}=\pm \sin 2\theta $). So we have 
\begin{equation*}
V_{1}^{+}=u^{+}T+w_{1}^{+}B,V_{1}^{-}=u^{-}T+w_{1}^{-}B,V_{2}^{\pm }=T
\end{equation*}%
and \ 
\begin{eqnarray*}
\beta _{2}^{\pm } &=&\left( (\lambda -a)\cos \frac{s}{c},(\lambda -a)\sin 
\frac{s}{c},\frac{bs}{c}\right)  \\
\beta _{1}^{\pm } &=&\left( \left( w_{2}^{\pm }b-u_{2}a-\lambda \right) \cos 
\frac{s}{c},\left( w_{2}^{\pm }b-u_{2}a-\lambda \right) \sin \frac{s}{c},%
\frac{w_{2}^{\pm }a+u_{2}b}{c}s\right) .
\end{eqnarray*}%
$\left( \alpha \right) $ is $\tan \theta -$Bertrand curve and the curves $%
\beta _{1}^{+},\beta _{1}^{-},\beta _{2}^{+},\beta _{2}^{-}$ are all $\tan
\theta -$Bertrand curve mate of $\left( \alpha \right) $.
\end{example}
\end{definition}

Let $\left( \gamma \right) $ be \ Bertrand curve and $\left\{ T,N,B,\kappa
,\ \tau \right\} $ be Serret-Frenet apparatus of the curves. We suppose that 
$v=0$, then we can see that $\lambda $ is non-zero constant$.$ We can obtain 
\begin{equation*}
V(s)=uT(s)+wB(s)
\end{equation*}%
where $u^{2}+w^{2}=1$ and $u,w\in R$. Integral curve of $V$ is defined as%
\begin{equation}
\gamma _{V}(s)=\int V(s)ds  \label{a}
\end{equation}%
Let $\left\{ T_{V},N_{V},B_{V},\kappa _{V},\ \tau _{V}\right\} $ \ be
Serret-Frenet apparatus of the curves $\left( \gamma _{V}\right) $. From
equation (\ref{a}), we have 
\begin{equation}
T_{V}=V(s)=uT(s)+wB(s)  \label{b}
\end{equation}%
From equation (\ref{b}), we can see that the curves $\left( \gamma \right) $
and $\left( \gamma _{V}\right) $ have a same arc-length parameter. From
equation (\ref{b}), we have 
\begin{equation*}
\kappa _{V}N_{V}=\left( u\kappa -w\tau \right) N
\end{equation*}%
So we obtain%
\begin{equation}
\ N_{V}=N  \label{c}
\end{equation}%
and 
\begin{equation}
\kappa _{V}=u\kappa -w\tau   \label{d}
\end{equation}%
From equation (\ref{b}) and (\ref{c}), we get 
\begin{equation}
B_{V}=T_{V}\times N_{V}=-wT(s)+uB(s)  \label{e}
\end{equation}%
From equation (\ref{e}), we have 
\begin{equation}
-\ \tau _{V}\ N_{V}=-\left( w\kappa +u\tau \right) N  \label{f}
\end{equation}%
So we obtain%
\begin{equation}
\tau _{V}\ =w\kappa +u\tau   \label{g}
\end{equation}%
From equation (\ref{d}) and (\ref{g}), we get%
\begin{equation}
\kappa \ =u\kappa _{V}+w\tau _{V}  \label{m}
\end{equation}%
and 
\begin{equation}
\tau =-w\kappa _{V}+u\tau _{V}  \label{n}
\end{equation}%
In this case, we can obtain following theorem.

\begin{theorem}
$\left( \gamma \right) $ is the Bertrand (resp. $B$-Bertrand ) curve if and
only if $\left( \gamma _{V}\right) $ is the Bertrand (resp. $B$-Bertrand )
curve.

\begin{proof}
We suppose that $\left( \gamma \right) $ is the Bertrand curve. There is $\
\lambda ,\mu \in R$ such that 
\begin{equation*}
\lambda \kappa +\mu \tau =u-w\cot \theta =1
\end{equation*}%
where $\lambda \neq 0$ , $\mu =\lambda \cot \theta $ \ and $\ \theta $\ is
constant angle between $T$ and $\overline{T}$. From equation (\ref{d}) and (%
\ref{g}), we get%
\begin{equation}
\lambda \kappa +\mu \tau =\overline{\lambda }\kappa _{V}+\overline{\mu }\tau
_{V}=1  \label{h}
\end{equation}%
where 
\begin{equation}
\overline{\lambda }=\lambda u-\mu w=\lambda \neq 0  \label{k}
\end{equation}%
and 
\begin{equation}
\overline{\mu }=\lambda w+\mu u=\lambda \cot (-\theta )  \label{l}
\end{equation}%
So, $\left( \gamma _{V}\right) $ is the Bertrand curve. Similarly, other
cases can be proved by the same method.
\end{proof}
\end{theorem}

Let $\left( \alpha \right) $ be \ Bertrand curve and $\left\{ T,N,B,\kappa
,\ \tau \right\} $ \ is Serret-Frenet apparatus of the curves. There is $%
\lambda $, $\mu \in 
%TCIMACRO{\U{211d} }%
%BeginExpansion
\mathbb{R}
%EndExpansion
$ \ such that $\lambda \kappa +\mu \tau =1$ where $\mu =\lambda \tan \theta $%
. If $\overline{\lambda }=\lambda a$ and $\overline{\mu }=\mu a$, then $%
\overline{\lambda }\kappa +\overline{\mu }\tau =a$. So we have%
\begin{equation}
u^{\pm }(a)=\frac{a\tan \theta \pm \sqrt{1+\left( \tan \theta \right)
^{2}-a^{2}}}{1+\left( \tan \theta \right) ^{2}}  \label{r}
\end{equation}%
and \ $w_{1}^{\pm }(a)=\pm \sqrt{1-\left( u^{+}(a)\right) ^{2}}$, $%
w_{2}^{\pm }(a)=\pm \sqrt{1-\left( u^{-}(a)\right) ^{2}}$. In this case, we
can give following definition.

\begin{definition}
Let $\left( \alpha \right) $ be \ Bertrand curve. So, we can define a
surface as 
\begin{eqnarray*}
\varphi _{1}^{+}(t,s) &=&\int V_{1}^{+}ds+\lambda N,\varphi
_{1}^{-}(t,s)=\int V_{1}^{-}ds+\lambda N \\
\varphi _{2}^{+}(t,s) &=&\int V_{2}^{+}ds+\lambda N,\varphi
_{2}^{-}(t,s)=\int V_{2}^{-}ds+\lambda N
\end{eqnarray*}%
where 
\begin{eqnarray*}
V_{1}^{+}(t,s)
&=&u^{+}(t)T(s)+w_{1}^{+}(t)B(s),V_{1}^{-}(t,s)=u^{+}(t)T(s)+w_{1}^{-}(t)B(s),
\\
V_{2}^{+}(t,s)
&=&u^{-}(t)T(s)+w_{2}^{+}(t)B(s),V_{2}^{-}(t,s)=u^{-}(t)T(s)+w_{2}^{-}(t)B(s),
\end{eqnarray*}%
$u^{\pm }(t)=\dfrac{t\tan \theta \pm \sqrt{1+\left( \tan \theta \right)
^{2}-t^{2}}}{1+\left( \tan \theta \right) ^{2}}$, \ $w_{1}^{\pm }(t)=\pm 
\sqrt{1-\left( u^{+}(t)\right) ^{2}}$ and $w_{2}^{\pm }(t)=\pm \sqrt{%
1-\left( u^{-}(t)\right) ^{2}}$. These surfaces are said Bertrand surface
taken by $\left( \alpha \right) $.
\end{definition}

\begin{example}
\ If \ $\alpha (s)=(\cos \frac{s}{\sqrt{2}},\sin \frac{s}{\sqrt{2}},\frac{s}{%
\sqrt{2}})$, then we have 
\begin{eqnarray*}
T &=&(-\frac{1}{\sqrt{2}}\sin \frac{s}{\sqrt{2}},\frac{1}{\sqrt{2}}\cos 
\frac{s}{\sqrt{2}},\frac{1}{\sqrt{2}}) \\
N &=&(-\cos \frac{s}{\sqrt{2}},-\sin \frac{s}{\sqrt{2}},0) \\
B &=&(\frac{1}{\sqrt{2}}\sin \frac{s}{\sqrt{2}},-\frac{1}{\sqrt{2}}\cos 
\frac{s}{\sqrt{2}},\frac{1}{\sqrt{2}})
\end{eqnarray*}%
and $\kappa =\frac{1}{2}$and $\tau =\frac{1}{2}$. If $\ \overline{\lambda }%
=\lambda \tan \theta =1$ and $\overline{\mu }=\lambda =1$, then we have $%
\overline{\lambda }\kappa +\overline{\mu }\tau =1=f$ and $\tan \theta =1$.
From (\ref{r}), we have 
\begin{equation*}
u(t)=\frac{t+\sqrt{2-t^{2}}}{2}
\end{equation*}%
and 
\begin{equation*}
w(t)=\frac{1}{\sqrt{2}}\sqrt{1-t\sqrt{2-t^{2}}}
\end{equation*}%
So, we have a surface $K\left( \alpha \right) $ as 
\begin{eqnarray*}
\varphi (t,s) &=&u(t)\int \ Tds+w(t)\int \ Bds+\lambda N \\
&=&\left( a(t)\cos \frac{s}{\sqrt{2}},a(t)\sin \frac{s}{\sqrt{2}},b(t)\frac{s%
}{\sqrt{2}}\right) 
\end{eqnarray*}%
where \ $a(t)=u(t)-w(t)-1$ and $b(t)=u(t)+w(t)$. The surface $K\left( \alpha
\right) $ is Bertrand surface taken by $\left( \alpha \right) $.
\end{example}

\begin{definition}
Let $\left( \left( \alpha \right) ,\left( \beta \right) \right) $ be \
Bertrand curve mate and $K\left( \alpha \right) $ and $K\left( \beta \right) 
$ be Bertrand surface. $\left( K\left( \alpha \right) ,K\left( \beta \right)
\right) $ is said Bertrand surface mate taken by $\left( \left( \alpha
\right) ,\left( \beta \right) \right) $ Bertrand mate.
\end{definition}

We can define a set as 
\begin{equation*}
K=\left\{ K\left( \alpha \right) /\left( \alpha \right) \text{ is a Bertrand
curve}\right\} \text{.}
\end{equation*}%
So we can define a equivalence relation on $K$ such that $K\left( \alpha
\right) \sim K\left( \beta \right) $ if and only if $\left( \left( \alpha
\right) ,\left( \beta \right) \right) $ is a\ Bertrand mate.

In Example 3.3, if $\left( \beta \right) $ is Bertrand mate of $\left(
\alpha \right) $, then we have 
\begin{eqnarray*}
\beta (s) &=&\alpha (s)+\lambda N \\
&=&(\cos \frac{s}{\sqrt{2}},\sin \frac{s}{\sqrt{2}},\frac{s}{\sqrt{2}}%
)+(-\cos \frac{s}{\sqrt{2}},-\sin \frac{s}{\sqrt{2}},0) \\
&=&(0,0,\frac{s}{\sqrt{2}})
\end{eqnarray*}

\begin{example}
If \ $\beta (s)=(0,0,\frac{s}{\sqrt{2}})$, then we have 
\begin{eqnarray*}
T &=&(0,0,1) \\
N &=&(-\cos \frac{s}{r},-\sin \frac{s}{r},0) \\
B &=&(\sin \frac{s}{r},-\cos \frac{s}{r},0)
\end{eqnarray*}%
and $\kappa =0$ and $\tau =r$. If $\ \overline{\lambda }=\lambda \tan \theta
=0$ and $\overline{\mu }=\lambda =r$, then we have $\overline{\lambda }%
\kappa +\overline{\mu }\tau =1=f$ and $\tan \theta =0$. From (\ref{r}), we
have $u(t)=\sqrt{1-t^{2}}$ and \ $w(t)=t$. So, we have a surface $K\left(
\alpha \right) $ as 
\begin{eqnarray*}
\varphi (t,s) &=&u(t)\int \ Tds+w(t)\int \ Bds+\lambda N \\
&=&\left( -\left( r+\frac{t}{r}\right) \cos \frac{s}{\sqrt{2}},-\left( r+%
\frac{t}{r}\right) \sin \frac{s}{\sqrt{2}},\sqrt{1-t^{2}}s\right) 
\end{eqnarray*}%
The surface $K\left( \beta \right) $ is Bertrand surface taken by $\left(
\beta \right) $. Furhermore, $\left( K\left( \alpha \right) ,K\left( \beta
\right) \right) $ is Bertrand surface mate taken by $\left( \left( \alpha
\right) ,\left( \beta \right) \right) $.
\end{example}

\section{\protect\bigskip Construction of the Bertrand and the
Principal-donor curves}

In $3$-Euclidean spaces, let $M$ be a regular spherical curve with
coordinate neighborhood $(I,\gamma )$. So we can define a curve $K$ ($\alpha
:I\rightarrow E^{3}$) such that 
\begin{equation}
\alpha (s)=\int S_{M}(s)\gamma (s)ds
\end{equation}%
where $S_{M}:I\longrightarrow 
%TCIMACRO{\U{211d} }%
%BeginExpansion
\mathbb{R}
%EndExpansion
$ $\left( s\rightarrow S_{M}(s)\right) $ is differentiable function.

\begin{lemma}
(\cite{CC}) The curve $K$ is spherical curve if and only if 
\begin{equation*}
S_{M}(s)=\left\Vert \gamma ^{\prime }(s)\right\Vert \cos \left( \overset{s}{%
\underset{0}{\dint }}\frac{\det (\gamma (s),\gamma ^{\prime }(s),\gamma
^{\prime \prime }(s))}{\left\Vert \gamma ^{\prime }(s)\right\Vert ^{2}}%
ds+\theta _{0}\right) .
\end{equation*}
\end{lemma}

Let $\ M$ \ be a unit-speed curve with coordinate neighborhood $(I,\gamma )$
and $\left\{ T,N,B,\kappa ,\ \tau \right\} $ be Serret-Frenet apparatus of
the curves. Because of unit-speed curve, we have $\left\Vert \gamma ^{\prime
}(s)\right\Vert =1$. From lemma 4.1, there is \ $S_{T}:I\longrightarrow 
%TCIMACRO{\U{211d} }%
%BeginExpansion
\mathbb{R}
%EndExpansion
$ differentiable function such that 
\begin{equation*}
\left\Vert \int S_{T}(s)\gamma ^{\prime }(s)ds\right\Vert =1
\end{equation*}%
where 
\begin{equation}
S_{T}(s)=\kappa (s)\cos \left( \overset{s}{\underset{0}{\dint }}\tau
(u)du+\theta _{0}\right)  \label{4.8}
\end{equation}%
If we define a curve $K$ with coordinate neighborhood $(I,\beta )$ such that 
\begin{equation*}
\beta ^{\prime }(s)=\int S_{T}(s)\gamma ^{\prime }(s)ds
\end{equation*}%
then we have 
\begin{equation}
\beta ^{\prime \prime }(s)=S_{T}(s)\gamma ^{\prime }(s)  \label{4.9}
\end{equation}%
Arc-lenght parameter of $M$ and $K$ are same. Let $\left\{ \overline{T},%
\overline{N},\overline{B},\overline{\kappa },\overline{\tau }\right\} $ be
Serret-Frenet apparatus of the curves where 
\begin{equation}
S_{T}(s)=\overline{\kappa }(s)=\kappa (s)\cos \left( \overset{s}{\underset{0}%
{\dint }}\tau (u)du+\theta _{0}\right)  \label{4.91}
\end{equation}%
and 
\begin{equation}
\overline{\tau }(s)=\kappa (s)\sin \left( \overset{s}{\underset{0}{\dint }}%
\tau (u)du+\theta _{0}\right) \text{.}  \label{4.10}
\end{equation}%
From (\ref{4.9}), we have 
\begin{equation}
\overline{N}(s)=\varepsilon T(s)  \label{4.11}
\end{equation}%
where $\varepsilon =\pm 1$. From equation (\ref{4.9}) and (\ref{4.11}), we
can see that principal normal of $K$ and tangent of $M$ is colinear.

\begin{theorem}
$K$ is principal-donor curve of $M$ if and only if $M$ is principal
direction curve of $K$.

\begin{proof}
From equation of $\beta ^{\prime \prime }(s)=S_{T}(s)\gamma ^{\prime }(s)$,
we have \ $\overline{N}=T$. So we can write%
\begin{equation}
\overline{T}=vN+wB  \label{1}
\end{equation}%
where $v^{2}+w^{2}=1$. If we take the derivative of equation (\ref{1}), we
have 
\begin{equation}
\overline{\kappa }\overline{N}=-\kappa vT+\left( v^{\prime }-w\tau \right)
N+\left( w^{\prime }+\tau v\right) B  \label{2}
\end{equation}%
So we have 
\begin{equation}
\overline{\kappa }=-\kappa v\text{, \ \ }v^{\prime }-w\tau =0\text{, \ }%
w^{\prime }+\tau v=0  \label{3}
\end{equation}%
From equation (\ref{3}), \ we have%
\begin{equation*}
v=-\cos \left( \overset{s}{\underset{0}{\dint }}\tau (u)du+\theta _{0}\right)
\end{equation*}%
and 
\begin{equation*}
w=\sin \left( \overset{s}{\underset{0}{\dint }}\tau (u)du+\theta _{0}\right)
\end{equation*}%
Thus, we have 
\begin{equation}
\overline{T}=-\cos \left( \overset{s}{\underset{0}{\dint }}\tau (u)du+\theta
_{0}\right) N+\sin \left( \overset{s}{\underset{0}{\dint }}\tau (u)du+\theta
_{0}\right) B  \label{4}
\end{equation}%
From equation (\ref{4}), $K$ is principal-donor curve of $M$. Conversely,
from equation $\beta ^{\prime \prime }(s)=\overline{\kappa }\gamma ^{\prime
}(s)$, we have $\overline{N}=T$. So we have $\gamma (s)=\overset{s}{\underset%
{0}{\dint }}\overline{N}(u)du$.
\end{proof}
\end{theorem}

\begin{theorem}
In 3-Euclidean spaces, let $M$, $K$ be a unit-speed curves with unit
coordinate neighborhoods $(I,\gamma )$ and $(I,\beta )$, respectively, such
that 
\begin{equation*}
\beta ^{\prime \prime }(s)=S_{T}(s)\gamma ^{\prime }(s)
\end{equation*}%
where%
\begin{equation}
S_{T}(s)=\kappa (s)\cos \left( \overset{s}{\underset{0}{\dint }}\tau
(u)du+\theta _{0}\right) \text{.}  \label{4.14}
\end{equation}
Following equations are equivalent,

1) $M$ is spherical curve

2) $K$ is Bertrand curve or $B$-Bertrand curve.

\begin{proof}
If $M$ is spherical curve, from (\ref{1.4}) we have 
\begin{equation*}
\frac{1}{\kappa (s)}=A\cos \left( \overset{s}{\underset{0}{\dint }}\tau
(u)du+\theta _{0}\right) 
\end{equation*}%
where $A$ and $\theta _{0}$ are constant. If $\ \theta _{0}=0$, then we have 
$\overline{\kappa }(s)=A(cont.)$. So $K$ is Salkowski curve. If $\ \theta
_{0}=-\frac{\pi }{2}$, then we have $\overline{\tau }(s)=A(const.)$. So, $K$
is anti-Salkowski curve. If $\ \theta _{0}\neq 0$ and $\theta _{0}\neq -%
\frac{\pi }{2}$, then we have 
\begin{equation*}
1=A\cos \left( \theta _{0}\right) \epsilon \kappa (s)\cos \theta (s)-A\sin
\left( \theta _{0}-c\right) \epsilon \kappa (s)\sin \theta (s)
\end{equation*}%
where $\theta (s)=\overset{s}{\underset{0}{\dint }}\tau (u)du$. So we have $%
\lambda \overline{\kappa }(s)+\mu \overline{\tau }(s)=1$ where 
\begin{equation*}
\lambda =A\cos \left( \theta _{0}\right) =const.\text{ and }\mu =-A\sin
\left( \theta _{0}\right) =const.
\end{equation*}%
Furthermore, $K$ is proper BertranD curve. Conversely, if $K$ is Bertrand
curve or $B$-Bertrand curve, there are $\lambda $,$\mu \in 
%TCIMACRO{\U{211d} }%
%BeginExpansion
\mathbb{R}
%EndExpansion
$ such that $\lambda \overline{\kappa }(s)+\mu \overline{\tau }(s)=1$. So we
have 
\begin{equation*}
1=\lambda \kappa (s)\cos \left( \overset{s}{\underset{0}{\dint }}\tau
(u)du+\theta _{0}\right) +\mu \kappa (s)\sin \left( \overset{s}{\underset{0}{%
\dint }}\tau (u)du+\theta _{0}\right) 
\end{equation*}%
and 
\begin{equation*}
\frac{1}{\kappa (s)}=A\left( \cos \left( \theta _{0}\right) \cos \left( 
\overset{s}{\underset{0}{\dint }}\tau (u)du\right) +\sin \left( \theta
_{0}\right) \sin \left( \overset{s}{\underset{0}{\dint }}\tau (u)du\right)
\right) 
\end{equation*}%
where $\theta _{0}\in \lbrack 0,2\pi ]$ , cos$\theta _{0}=\frac{\lambda }{%
\sqrt{\lambda ^{2}+\mu ^{2}}}$, sin$\theta _{0}=-\frac{\mu }{\sqrt{\lambda
^{2}+\mu ^{2}}}$ and $A=\sqrt{\lambda ^{2}+\mu ^{2}}$.\ So we have $\frac{1}{%
\kappa (s)}=A\cos \left( \overset{s}{\underset{0}{\dint }}\tau (u)du+\theta
_{0}\right) $.
\end{proof}
\end{theorem}

\begin{example}
In 3-Euclidean spaces, let $M$ be a circle where $\gamma (s)=(r\cos
ws,-r\sin ws,0)$ and $w=\frac{1}{r}$. So we have $\kappa =w$ and $\tau =0$.
From (\ref{4.91}) and (\ref{4.10}), we obtain $\overline{\kappa }%
(s)=\epsilon w\cos c_{0}=A(const.)$ and $\overline{\tau }(s)=\overline{%
\epsilon }w\sin c_{0}=B(const.)$. So, from (\ref{4.9}) , we get 
\begin{equation*}
\beta ^{\prime \prime }(s)=A(-\sin ws,-\cos ws,0)
\end{equation*}%
If we integrate of the above equation, we have a curve $K$ as 
\begin{equation*}
\beta (s)=(Ar^{2}\cos \left( \epsilon ws\right) ,Ar^{2}\sin \left( \epsilon
ws\right) ,bws+c_{1})
\end{equation*}%
Because $K$ is unit spreed curve, we have $b=r\sin c_{0}$. If $%
a=Ar^{2}=\epsilon r\cos c_{0}$ and $c_{1}=0$, we have%
\begin{equation*}
\beta (s)=(a\cos \left( ws\right) ,a\sin \left( ws\right) ,bws)
\end{equation*}%
where $r=\sqrt{a^{2}+b^{2}}$.
\end{example}

\end{document}